\documentclass{article}

\newtheorem{fed}{\textbf{Definition}}[section]
\newtheorem{thm}[fed]{\textbf{Theorem}}
\newtheorem{lemma}[fed]{\textbf{Lemma}}

\newtheorem{prop}[fed]{\textbf{Proposition}}
\newtheorem{cor}[fed]{\textbf{Corollary}}
\usepackage{amssymb,bbm,graphicx,epsfig,psfrag,epic,eepic,latexsym}
\usepackage{amsmath}
\usepackage{mathrsfs}
\usepackage[dvips]{color}

\begin{document}
\title{State dependent Hamiltonian delay equations and Neumann one-forms}
\author{Urs Frauenfelder}
\maketitle
\begin{abstract}
In this note we study critical points of a variation of the action functional of classical mechanics, where the 
Hamiltonian term is retarded. Following a more than hundert and fifty year old paper by Carl Neumann we as well 
introduce Taylor approximations to this functional in terms of the fine structure constant. We see how in first 
order in the fine structure constant these functionals are related to time dependent perturbations of the symplectic
form. This observation allows us to prove Arnold type conjectures in first order in the fine structure constant. 
\end{abstract}

\tableofcontents

\section{Introduction}

In this note we explore action functionals which admit a state dependent delay. For exact symplectic manifolds $(M,\omega=d\lambda)$ such action functionals have the following form. If $S^1=\mathbb{R}/\mathbb{Z}$ is the circle,
$H \in C^\infty(M \times S^1,\mathbb{R})$ is a Hamiltonian depending periodically on time, $F \in C^\infty(M,\mathbb{R})$
is another smooth function, referred to as the \emph{retardation} and $\mathcal{L}=C^\infty(S^1,M)$ is the free loop space,
then the action functional we are interested in is
$$\mathcal{A} \colon \mathcal{L} \to \mathbb{R},\quad v \mapsto \int_{S^1} v^*\lambda-\int_0^1 H\Big(v\big(t+F(v(t))\big),t+F(v(t))\Big)dt.$$
If $F=0$ this is just the action functional of classical mechanics. The study of action functionals of this kind is motivated by an old paper of Carl Neumann \cite{neumann1} from 1868. In this paper Neumann studied an action functional in which the Coulomb potential is evaluated at the retarded time 
$$\tau=t-\frac{|q(t)|}{c}.$$ 
The author is currently exploring together with Peter Albers and Felix Schlenk the possibility to construct nonlocal Floer homologies and apply them to Hamiltonian delay equations, see \cite{albers-frauenfelder-schlenk1, albers-frauenfelder-schlenk2, albers-frauenfelder-schlenk3}. He expects that action functionals of the above form give rise to a rich source of examples on which one can check how far nonlocal Floer homology can be developed but also where its limits lie. For that purpose we formulate some Post-Arnold questions motivated by the famous Arnold conjecture about topological lower bounds on the number of periodic orbits of Hamiltonian systems. 
\\ \\
Following an idea of Neumann we consider a whole one-parameter family of action functionals, namely for $\alpha \in [0,\infty)$ we look at
$$\mathcal{A}_\alpha \colon \mathcal{L} \to \mathbb{R},\quad v \mapsto \int_{S^1} v^*\lambda-\int_0^1 H\Big(v\big(t+\alpha F(v(t))\big),t+\alpha F(v(t))\Big)dt.$$
For $\alpha=0$ we recover the action functional of classical mechanics. In the above interpretation of retarded time the parameter $\alpha$ corresponds to $\tfrac{1}{c}$. Interpreting $c$ as the speed of light its inverse in atomic units corresponds to the fine structure constant, whence we refer to $\alpha$ as the fine structure constant. As in the paper by Neumann \cite{neumann1} we consider as well the Taylor polynomial of $\mathcal{A}_\alpha$ in $\alpha$ at the origin $\alpha=0$. In the case of the Coulomb potential Neumann discovered that by considering the quadratic Taylor polynomial one discovers Weber's electrodynamics, see also \cite{frauenfelder-weber}.
\\ \\
In the general Hamiltonian case considered in this paper, already the linear Taylor polynomial shows interesting phenomena. It is related to time-dependent perturbations of the symplectic form. This interpretation allows us to prove for small enough fine structure constant an Arnold type conjecture in first order in the fine structure constant. The precise statement of this result is Theorem~\ref{main} below. Because we are just working in the Taylor approximation, the corresponding Floer homology is still local, however, in contrast to the classical case it requires time-dependent symplectic forms. This in particular requires some careful estimates on the energy of gradient flow lines in order to construct the Floer homology. 
\\ \\
In the appendix we study retardations which are not local anymore but might depend on the whole loop. We see how to first order in the fine structure constant this is related to nonlocal symplectic forms on the free loop space. 
\\ \\
\emph{Acknowledgement:} The author would like to thank Felix Schlenk for useful comments. He gratefully acknowledges support from the grant FR 2637/2-1 funded by the German Research Foundation (DFG).

\section{Neumann one-forms and its critical points}

\subsection{Periodic orbits of a Hamiltonian vector field}

In this section we recall how periodic orbits of a time-periodic Hamiltonian vector field can be obtained variationally as critical points of a closed one-form on the free loop space.
\\ \\
Assume that $(M,\omega)$ is a symplectic manifold. We abbreviate by 
$$\mathcal{L}:=C^\infty(S^1,M)$$
the free loop space of $M$ where $S^1=\mathbb{R}/\mathbb{Z}$ is the circle. On $\mathcal{L}$ we have the closed one-form
$$\mathfrak{a}_\omega \in \Omega^1(\mathcal{L})$$
which is defined as follows. If $v \in \mathcal{L}$ and $\hat{v} \in T_v \mathcal{L}=\Gamma(v^*TM)$ is a smooth vector field along $v$, then
\begin{equation}\label{oneform}
\mathfrak{a}_\omega(v) \hat{v}:=\int_{S^1} v^* \iota_{\hat{v}}\omega=\int_0^1 \omega_{v(t)}\big(\hat{v}(t),\partial_t v(t)\big)dt.
\end{equation}
Suppose that
$$H \in C^\infty(M\times S^1,\mathbb{R})$$
is a smooth function on $M \times S^1$. For $t \in S^1$ abbreviate
$$H_t:=H(\cdot,t) \in C^\infty(M,\mathbb{R}).$$
We refer to $H$ as the Hamiltonian, which periodically depends on time. Its
mean value
\begin{equation}\label{mean}
\mathcal{H} \colon \mathcal{L}\to \mathbb{R},\quad v \mapsto \int_0^1 H(v(t),t)dt
\end{equation}
is a smooth function on the free loop space. The critical points of the closed one form
$$\mathfrak{a}_\omega-d\mathcal{H} \in \Omega^1(\mathcal{L})$$
correspond to periodic orbits of period one of the time dependent Hamiltonian vector field $X_{H_t}$ of $H_t$ implicitly defined by the condition
$$dH_t=\omega(\cdot, X_{H_t}).$$
Indeed, suppose that
$$v \in \mathrm{crit}\big(\mathfrak{a}_\omega-d\mathcal{H}\big).$$
Then for every $\hat{v} \in T_v \mathcal{L}$
$$0=(\mathfrak{a}_\omega-d\mathcal{H}\big)(v)\hat{v}=\int_0^1 \omega(\hat{v},\partial_t v)dt-\int_0^1
dH_t(v)\hat{v}dt=\int_0^1 \omega\big(\hat{v},\partial_t v-X_{H_t}(v)\big)dt$$
implying that $v$ is a solution of the ODE
\begin{equation}\label{periodic}
\partial_t v=X_{H_t}(v)
\end{equation}
i.e., $v$ is a periodic orbit of $X_{H_t}$ of period one. From Floer homology we know that on closed symplectic manifolds
periodic orbits always exist.

\subsection{Neumann one-forms}

In his paper from 1868 \cite{neumann1} Carl Neumann considered action functionals which admit a state dependent delay. 
In his example he considered a retarded Coulomb potential and showed how Weber's force law can be deduced in this way, see also \cite{frauenfelder-weber}. A modern interpretation is that the proton does not attract the electron
instantaneously but that the force is transmitted by virtual photons of finite speed, which leads to a retarded potential. 
\\ \\
We next explain how to incorporate such a retardation in a Hamiltonian set-up. 
Suppose that 
$$F \in C^\infty(M,\mathbb{R})$$
is another smooth function on $M$ which we refer to as a \emph{retardation}. Let
$$\alpha \in [0,\infty)$$ 
be a parameter which we refer to as the \emph{fine structure constant}.
Consider the perturbation
$$\mathcal{H}_\alpha \colon \mathcal{L} \to \mathbb{R}$$
of $\mathcal{H}$ which for $v \in \mathcal{L}$ is given by
$$\mathcal{H}_\alpha(v):=\int_0^1 H\Big(v\big(t+\alpha F(v(t))\big),t+\alpha F(v(t))\Big)dt.$$
Here is an example which as well explains the terminology. Suppose that $M=T^*\mathbb{R}^n=\mathbb{R}^n\times\mathbb{R}^n$
and let $(q,p) \in \mathbb{R}^n \times \mathbb{R}^n$ where we think of $q$ as the position and of $p$ as the momentum. 
Let
$$F(q,p)=-|q|.$$
In atomic units the fine structure constant $\alpha$ corresponds to the inverse of the speed of light
$$\alpha=\frac{1}{c}$$
so that $H \circ v$ is evaluated at the retarded time
$$\tau=t-\frac{|q|}{c}$$
where one can imagine that the retardation is due to the fact that a virtual photon has to travel from the origin to
$q$ with speed $c$.
Note that the case where the fine structure is zero corresponds to the case where the virtual photon has infinite speed so that indeed we have 
$$\mathcal{H}_0=\mathcal{H}.$$ 
We refer to 
$$\mathfrak{a}_\omega-d\mathcal{H}_\alpha \in \Omega^1(\mathcal{L})$$
as a Neumann one-form on the free loop space. 

\subsection{The autonomous case}

In this section we discuss the case where $H$ is independent of $t \in S^1$ and therefore just a smooth function on $M$. 
Under this assumption the formula for $\mathcal{H}_\alpha(v)$ simplifies to
$$\mathcal{H}_\alpha(v):=\int_0^1 H\Big(v\big(t+\alpha F(v(t))\big)\Big)dt.$$
Note that the circle $S^1$ acts on the free loop space $\mathcal{L}$ by rotating the domain, i.e., if $r \in S^1$ and
$v \in \mathcal{L}$, then
$$r_* v(t)=v(t+r), \quad t \in S^1.$$
In the autonomous case, the function $\mathcal{H}_\alpha$ is invariant under the circle action on the free loop space, i.e.
$$\mathcal{H}_\alpha(r_* v)=\mathcal{H}_\alpha(v), \quad r \in S^1,\,\,v\in \mathcal{L}.$$
Because the one-form $\mathfrak{a}_\omega$ is invariant as well we obtain the following lemma.
\begin{lemma}
In the autonomous case the critical set $\mathrm{crit}(\mathfrak{a}_\omega-d\mathcal{H}_\alpha)$ is invariant under
the circle action on the free loop space. In particular, critical points either come in circle families or correspond to
fixed points of the circle action, i.e., constant loops.
\end{lemma}
To compute the differential of $\mathcal{H}_\alpha$ at $v \in \mathcal{L}$ we pick $\hat{v} \in T_v\mathcal{L}$ and
compute 
\begin{eqnarray}\label{diff1}
d\mathcal{H}_\alpha(v)\hat{v}&=&\int_0^1 dH\Big(v\big(t+\alpha F(v(t))\big)\Big)\hat{v}\big(t+\alpha F(v(t))\big)dt\\ \nonumber
& &+
\int_0^1 dH\Big(v\big(t+\alpha F(v(t))\big)\Big)\partial_t v\big(t+\alpha F(v(t))\big)\alpha dF(v(t))\hat{v}(t)dt\\ \nonumber
&=&\int_0^1 \omega\bigg(\hat{v}\big(t+\alpha F(v(t))\big), X_H\Big(v\big(t+\alpha F(v(t))\big)\Big)\bigg)dt\\ \nonumber
& &+
\int_0^1 \alpha dH\Big(v\big(t+\alpha F(v(t))\big)\Big)\partial_t v\big(t+\alpha F(v(t))\big)\omega\big(\hat{v}(t),X_F(v(t))\big)dt
\end{eqnarray}
\begin{fed}
Suppose that $v \in \mathcal{L}$. A time
$t_0 \in S^1$ is called \emph{lightlike} (with respect to $v$) if
\begin{equation}\label{lightlike}
\alpha dF(v(t_0))\partial_t v(t_0)=-1.
\end{equation}
\end{fed}
The reason for this terminology is the following. We can rewrite equation (\ref{lightlike}) as
$$-dF(v(t_0))\partial_t v(t_0)=\frac{1}{\alpha}$$
and by interpreting the inverse value of the fine structure constant as the velocity of light, this means that roughly
the velocity of $v$ at time $t_0$ corresponds to the velocity of light. 
\begin{fed}
A loop $v \in \mathcal{L}$ is called \emph{spacelike},
if
$$\alpha dF(v(t)) \partial_t v(t)>-1,\quad \forall\,\,t \in S^1.$$
\end{fed}
Note that for a spacelike loop no lightlike times exist. 
\\ \\
Given $v \in \mathcal{L}$, define 
$$\tau_v \colon S^1 \to S^1,\quad t \mapsto t+\alpha F(v(t))\,\,\mathrm{mod}\,\,1.$$
\begin{lemma}\label{sp}
Suppose that $v \in \mathcal{L}$ is spacelike. Then $\tau_v$ is a diffeomorphism of the circle. 
\end{lemma}
\textbf{Proof:} First observe that $\tau_v \colon S^1 \to S^1$ has always degree one, whether or not $v$ is spacelike. If $v$ is spacelike then
$$\frac{d\tau_v(t)}{dt}=1+\alpha dF(v(t))\partial_t v>0$$
so that $\tau_v$ is strictly monotone increasing and therefore a diffeomorphism. \hfill $\square$
\\ \\
In view of Lemma~\ref{sp} for spacelike $v \in \mathcal{L}$ we can rewrite (\ref{diff1})
\begin{eqnarray*}
d\mathcal{H}_\alpha(v)\hat{v}
&=&\int_0^1 \frac{d\tau_v^{-1}(t)}{dt}\omega\big(\hat{v}(t), X_H(v(t))\big)dt\\ 
& &+
\int_0^1 \alpha dH\Big(v\big(t+\alpha F(v(t))\big)\Big)\partial_t v\big(t+\alpha F(v(t))\big)\omega\big(\hat{v}(t),X_F(v(t))\big)dt.
\end{eqnarray*}
In view of this formula we obtain the following proposition.
\begin{prop}\label{del1}
A loop $v$ is a spacelike critical point of $\mathfrak{a}_\omega-d\mathcal{H}_\alpha$, if and only if $v$ is a solution
of the following delay equation
\begin{equation}\label{delay}
\partial_t v(t)=\frac{d\tau_v^{-1}(t)}{dt}X_H(v(t))+\alpha dH\Big(v\big(t+\alpha F(v(t))\big)\Big)\partial_t v\big(t+\alpha F(v(t))\big)X_F(v(t)).
\end{equation}
\end{prop}
An immediate consequence of this formula is
\begin{cor}
Each critical point of $H$ interpreted as a constant loop is a spacelike critical point of $\mathfrak{a}_\omega-d\mathcal{H}_\alpha$.
\end{cor}
This corollary has as a further consequence the following
\begin{cor}\label{autoarnold}
There are at least as many contractible spacelike critical points of $\mathfrak{a}_\omega-d\mathcal{H}_\alpha$
as the minimal number of critical points of a function on $M$.
\end{cor}

\subsection{Preservation of energy and the parallel case}

The results of this subsection are not used in the sequel and the hasty reader can skip it. To the leisurely reader,
however, the following question may come to mind: In the autonomous case, is there a pointwise preserved quantity for 
critical points of
$\mathfrak{a}_\omega-d\mathcal{H}_\alpha$? Indeed, for autonomous Hamiltonian systems the Hamiltonian is preserved along its flow and therefore also referred to as the energy. Starring at (\ref{delay})
we see that this continues to hold if $H$ and $F$ Poisson commute. Recall that the Poisson bracket on a symplectic manifold
$$\{\cdot,\cdot\} \colon C^\infty(M)\times C^\infty(M) \to C^\infty(M)$$
is defined as
$$\{H,F\}:=\omega(X_H,X_F)=dF(X_H)=-dH(X_F).$$
It is antisymmetric and therefore in particular,
$$\{H,H\}=0.$$
We say that $H$ and $F$ Poisson commute if $\{H,F\}=0$. 
\begin{lemma}
Suppose that $H$ and $F$ Poisson commute and that $v \in \mathcal{L}$ is a spacelike critical point of
$\mathfrak{a}_\omega-d\mathcal{H}_\alpha$. Then $H$ and $F$ are constant along $v$.
\end{lemma}
\textbf{Proof: } Given $v$ as in the lemma, we see from (\ref{delay}) that there exist smooth functions 
$$a,b \colon S^1 \to \mathbb{R}$$
such that 
$$\partial_t v(t)=a(t)X_H(v(t))+b(t)X_F(v(t)).$$
To show that $H$ is constant along $v$ we compute
\begin{eqnarray*}
\frac{d}{dt}H(v(t))&=&dH(v(t))\partial_t v(t)\\
&=&dH(v(t))\Big(a(t)X_H(v(t))+b(t)X_F(v(t))\Big)\\
&=&a(t)\{H,H\}(v(t))+b(t)\{F,H\}(v(t))\\
&=&0.
\end{eqnarray*}
Hence by the fundamental theorem of calculus $H \circ v$ is constant. That $F \circ v$ is constant follows in
the same way. \hfill $\square$
\\ \\
A special case where $F$ and $H$ Poisson commute is when the differential of $F$ is parallel to $H$, i.e., there exists
$G \in C^\infty(M,\mathbb{R})$ such that
\begin{equation}\label{parallel}
dF=GdH.
\end{equation}
In this case the Hamiltonian vector fields are parallel as well, namely
$$X_F=GX_H.$$
This happens for instance if $F$ can be written as
$$F=f \circ H$$
for some smooth function $f \colon \mathbb{R} \to \mathbb{R}$ so that, in particular, each level set of $H$ is contained in a level set of $F$. In this case (\ref{parallel}) holds for
$$G=f' \circ H.$$
In the parallel case equation (\ref{delay}) for spacelike critical points simplifies to
\begin{equation}\label{pardelay}
\partial_t v(t)=\bigg(\frac{d\tau_v^{-1}(t)}{dt}+\alpha dH\Big(v\big(t+\alpha F(v(t))\big)\Big)\partial_t v\big(t+\alpha F(v(t))\big)G(v(t))\bigg)X_H(v(t)).
\end{equation}
In the parallel case the velocity of $v$ is always parallel to the Hamiltonian vector field of $H$, but the factor
of proportionality can depend on the whole trajectory. In view of this property we have the following lemma.
\begin{lemma}
Suppose that $dF$ is parallel to $dH$, that $v$ is a spacelike critical point of $\mathfrak{a}_\omega-d\mathcal{H}_\alpha$
and that there exists $t_0 \in S^1$ such that $v(t_0) \in \mathrm{crit}(H)$. Then $v$ is constant.
\end{lemma}
\textbf{Proof:} By (\ref{pardelay}) we see that there exists a smooth function $f \colon S^1 \to \mathbb{R}$ such that
$v$ is a solution of the ODE 
\begin{equation}\label{ode}
\partial_t v(t)=f(t)X_H(v(t)).
\end{equation}
Because $v(t_0) \in \mathrm{crit}(H)$, we have $X_H(v(t_0))=0$, and so
the constant loop
$$v_0 \colon S^1 \to M, \quad t \mapsto v(t_0)$$
is a solution of (\ref{ode}) as well. Moreover, $v$ and $v_0$ coincide at $t_0$. By Lindel\"of's uniqueness theorem for
the initial value problem of an ODE we conclude that
$$v=v_0$$
so that in particular $v$ is constant. \hfill $\square$

\subsection{The nonautonomous case}

In the nonautonomous case the Hamiltonian $H$ is allowed to depend periodically on time. This has the effect that the partial derivative in time direction
$$\partial_t H \in C^\infty(M \times S^1,\mathbb{R})$$
gives rise to an additional term in (\ref{diff1}). Indeed, in the nonautonomous case the computation in (\ref{diff1}) has
to be replaced by
\begin{eqnarray}\label{diff2}
& &d\mathcal{H}_\alpha(v)\hat{v}\\ \nonumber
&=&\int_0^1 dH_{t+\alpha F(v(t))}\Big(v\big(t+\alpha F(v(t))\big)\Big)\hat{v}\big(t+\alpha F(v(t))\big)dt\\ \nonumber
& &+
\int_0^1 dH_{t+\alpha F(v(t))}\Big(v\big(t+\alpha F(v(t))\big)\Big)\partial_t v\big(t+\alpha F(v(t))\big)\alpha dF(v(t))\hat{v}(t)dt\\ \nonumber
& &+\int_0^1 \partial_t H\Big(v\big(t+\alpha F(v(t))\big),t+\alpha F(v(t))\Big)\alpha dF(v(t))\hat{v}(t)dt\\ \nonumber
&=&\int_0^1 \omega\bigg(\hat{v}\big(t+\alpha F(v(t))\big), X_{H_{t+\alpha F(v(t))}}\Big(v\big(t+\alpha F(v(t))\big)\Big)\bigg)dt\\ \nonumber
& &+
\int_0^1 \alpha dH_{t+\alpha F(v(t))}\Big(v\big(t+\alpha F(v(t))\big)\Big)\partial_t v\big(t+\alpha F(v(t))\big)\omega\big(\hat{v}(t),X_F(v(t))\big)dt\\ \nonumber
& &+\int_0^1 \alpha \partial_t H\Big(v\big(t+\alpha F(v(t))\big),t+\alpha F(v(t))\Big)\omega\big(\hat{v}(t),X_F(v(t))\big)
\end{eqnarray}
In view of this formula the generalization of Proposition~\ref{del1} to the nonautonomous case becomes
\begin{prop}\label{del2}
Suppose that $v$ is a spacelike critical point of $\mathfrak{a}_\omega-d\mathcal{H}_\alpha$, then $v$ is a solution
of the following delay equation
\begin{eqnarray}\label{delay2}
\partial_t v(t)&=&\frac{d\tau_v^{-1}(t)}{dt}X_{H_t}(v(t))\\ \nonumber
& &+\alpha dH_{t+\alpha F(v(t))}\Big(v\big(t+\alpha F(v(t))\big)\Big)\partial_t v\big(t+\alpha F(v(t))\big)X_F(v(t))
\\ \nonumber
& &+\alpha \partial_t H\Big(v\big(t+\alpha F(v(t))\big),t+\alpha F(v(t))\Big)X_F(v(t)).
\end{eqnarray}
\end{prop}

\section{Post-Arnold questions}\label{arnold}

In the 1960s, see also \cite[Appendix 9]{arnold}, Arnold formulated his famous conjecture that the number of periodic orbits of a Hamiltonian vector field on a closed symplectic manifold can be estimated from below by the minimal number of critical points of a function on this manifold. This conjecture was one of the important driving forces in the development of many new techniques in symplectic topology. It inspired Floer to the construction of his semi-infinite dimensional Morse homology, nowadays referred to as Floer homology \cite{floer1,floer2}, it led Fukaya and Ono to the discovery of Kuranishi structures \cite{fukaya-ono}, and it plays a motivating role in the development of the theory of polyfolds discovered by Hofer, Wysocki, and Zehnder \cite{filippenko-wehrheim, hofer-wysocki-zehnder}. To the author's knowledge the original version of the Arnold conjecture is just proved for symplectically aspherical manifolds, i.e., symplectic manifolds which have the property that the integral of the symplectic form over spheres vanishes. This result is due to Rudyak and Oprea
\cite{rudyak-oprea} and combines Floer theory with clever topological methods. For general closed symplectic manifolds the existence of one periodic orbit is nowadays always established and for generic Hamiltonians much stronger results hold true, namely topological lower bounds in terms of the sum of the Betti-numbers of the manifold, see for instance \cite[Chapter 12]{mcduff-salamon2}. Although the author is not aware of an explicit counterexample to the original version of the Arnold conjecture in the non symplectically aspherical case, it is doubtful if the original conjecture of Arnold holds true in this case. Indeed, as Floer homology teaches us in the non symplectically aspherical case, the topology of the manifold has to be replaced by its quantum topology, see for instance \cite{schwarz1}.
\\ \\
In the symplectically aspherical case, on the component of contractible loops $\mathcal{L}_0 \subset \mathcal{L}$ the 
one-form $\mathfrak{a}_\omega \in \Omega(\mathcal{L}_0)$ defined in (\ref{oneform}) is exact. Namely it arises as the differential of the area functional. Indeed, if $v \in \mathcal{L}_0$ we can find a filling disk
$$\overline{v} \colon D \to M$$
where $D:=\{z \in \mathbb{C}: |z| \leq 1\}$ is the closed unit ball in $\mathbb{C}$ such that
$$\overline{v}(e^{2 \pi i t})=v(t).$$
We then define the area functional
$$\mathcal{A}_0 \colon \mathcal{L}_0 \to \mathbb{R}, \quad v \mapsto \int_D \overline{v}^* \omega$$
and due to the assumption that our symplectic manifold is symplectically aspherical this functional is well defined independent of the choice of the filling disk. Indeed, if we choose another filling disk for $v$ we can glue the two filling disks along $v$
to obtain a sphere, and because the integral of the symplectic form over this sphere vanishes the integrals over each of the filling disks have to coincide. Note that
$$d \mathcal{A}_0=\mathfrak{a}_\omega$$
If $H \in C^\infty(M \times S^1,\mathbb{R})$ is a Hamiltonian depending periodically on time we set
$$\mathcal{A}_H \colon \mathcal{L}_0 \to \mathbb{R}, \quad \mathcal{A}_H:=\mathcal{A}_0-\mathcal{H}.$$
Therefore we can interpret contractible periodic orbits of the Hamiltonian vector field of $H$, i.e., solutions of (\ref{periodic}), as critical points of $\mathcal{A}_H$. Classical Floer homology is the semi-infinite dimensional Morse homology associated to $\mathcal{A}_H$ and therefore the periodic orbits that it detects have the additional property of being contractible. In the case that the symplectic manifold is not symplectically aspherical one has to consider a covering of the component of contractible loops on which the closed one-form $\mathfrak{a}_\omega$ becomes exact. Floer homology in this case is the semi-infinite dimensional Morse-Novikov homology associated to it. 
\\ \\
In the symplectically aspherical case we can associate to a periodic orbit $v$ its action $\mathcal{A}_H(v)$. A theorem due to Schwarz \cite{schwarz2} tells us that in the case where not through every point of $M$ passes a periodic orbit, there are at least two periodic orbits of different action. 
\\ \\
To state our Post-Arnold questions we assume that $(M,\omega)$ is a closed symplectically aspherical manifold.
We fix a Hamiltonian $H$ depending periodically on time and a retardation $F$. 
For the first question, if a positive result can be obtained in the symplectically aspherical case, then probably such a result can as well be extended to the more general case. However, the assertion probably has to be modified a bit taking the topology of the quantum cup product into account, and the methods require presumably heavy technology like Kuranishi structures or polyfolds.
\\ \\
\textbf{Question\,1:} \emph{Can the number of contractible spacelike critical points of the Neumann one form $\mathfrak{a}_\omega-d\mathcal{H}_\alpha$ for $\alpha \in [0,\infty)$ be estimated from below by $\mathrm{Crit}(M)$, the minimal number of critical points of a function on $M$?}
\\ \\
\textbf{Question\,2:} \emph{For $\alpha \in [0,\infty)$, in the case that not through every point of $M$ passes a spacelike contractible critical point of $\mathfrak{a}_\omega-d\mathcal{H}_\alpha$, do there always exist two contractible spacelike critical points of different action?}
\\ \\
Note that in the autonomous case Question\,1 has a positive answer due to Corollary~\ref{autoarnold}. Moreover, observe that
in the autonomous case the action of a critical point $v$ of $H$ interpreted as a constant loop is just
\begin{equation}\label{autact}
\mathcal{A}_H(v)=-H(v).
\end{equation}
If $H$ is constant, then all points of $M$ are critical. Hence to answer Question\,2 in the autonomous case we can assume without loss of generality that $H$ is not constant. But then if $v^+$ is a point at which $H$ attains its maximum and
$v^-$ is a point at which $H$ attains its minimum we have
$$H(v^+)>H(v^-)$$ 
and by (\ref{autact}) we found two critical points of $\mathcal{A}_H$ of different action. This shows that in the autonomous case also Question~2 admits a positive answer. 
\\ \\
Despite the fact that both questions admit a positive answer in the autonomous case, the author of this note is not too optimistic that this continuous to hold in the non autonomous case and suggests to a reader who really wants to know the answer to these questions to try first to construct a counterexample. We therefore weaken the two questions as follows. 
\\ \\
\textbf{Question\,$1^{-}$: } \emph{Does there exist a constant $\alpha_0=\alpha_0(H,F)>0$ such that for every $\alpha \in [0,\alpha_0]$
the number of contractible spacelike critical points of $\mathfrak{a}_\omega-d\mathcal{H}_\alpha$ is at least $\mathrm{Crit}(M)$?}
\\ \\
\textbf{Question\,$2^{-}$: } \emph{Does there exist a constant $\alpha_0=\alpha_0(H,F)>0$ with the property that for every $\alpha \in [0,\alpha_0]$ for which not through every point of $M$ passes a spacelike contractible critical point of
$\mathfrak{a}_\omega-d\mathcal{H}_\alpha$, do there always exist two contractible spacelike critical points of different action?}
\\ \\
We show in the next two sections that in first order of the fine structure  constant the answer to Question\,$1^{-}$ and Question\,$2^{-}$ is positive. For higher order in the fine structure constant we do not know.

\section{Taylor polynomials of Neumann one-forms}

\subsection{Taylor polynomials}

We fix $v \in \mathcal{L}$ and consider the Taylor polynomials at $\alpha=0$ for the function
$\alpha \mapsto \mathcal{H}_\alpha(v)$. To do that we abbreviate
\begin{equation}\label{hv}
H_v \in C^\infty([0,\infty)\times S^1,\mathbb{R}),\quad (\alpha,t) \mapsto H\Big(v\big(t+\alpha F(v(t))\big),
t+\alpha F(v(t))\Big).
\end{equation}
For $k \in \mathbb{N}_0$ we define
$$\mathcal{H}^k \colon \mathcal{L} \to \mathbb{R}$$
by
$$\mathcal{H}^k(v):=\frac{1}{k!}\int_0^1 \frac{\partial^k}{\partial \alpha^k}H_v(0,t)dt.$$ 
For $n \in \mathbb{N}_0$ we set
$$\mathcal{H}_{n,\alpha}:=\sum_{k=0}^n \alpha^k \mathcal{H}^k \colon \mathcal{L}\to \mathbb{R}.$$
Note that 
$$\mathcal{H}_{0,\alpha}=\mathcal{H}^0=\mathcal{H}.$$

\subsection{First order in the autonomous case and deformation of the symplectic form}

In this subsection we show that in first order in the fine structure constant the retardation can be incorporated in a deformation of the symplectic form.
\\ \\
Note that in the autonomous case (\ref{hv}) simplifies to
$$H_v \in C^\infty([0,\infty)\times S^1,\mathbb{R}),\quad (\alpha,t) \mapsto H\Big(v\big(t+\alpha F(v(t))\big)\Big).$$
In order to get the linear term in the fine structure constant we compute
$$\frac{\partial}{\partial \alpha} H_v(\alpha,t)=dH\Big(v\big(t+\alpha F(v(t))\big)\Big)\partial_t v
\big(t+\alpha F(v(t))\big)F(v(t)).$$
In particular, setting the fine structure constant equal to zero we obtain
$$\frac{\partial}{\partial \alpha} H_v(0,t)=F(v(t))dH(v(t))\partial_t v(t)=F(v(t))\frac{d}{dt}H(v(t)).$$
Hence if we introduce the one-form
$$\zeta=FdH \in \Omega^1(M)$$
we can write
$$\mathcal{H}^1(v)=\int_{S^1} v^* \zeta.$$
If $\hat{v} \in T_v \mathcal{L}$ we obtain by using Cartan's formula and Stokes' Theorem
\begin{eqnarray*}
d \mathcal{H}^1(v)\hat{v}&=&\int_{S^1}v^*L_{\hat{v}}\zeta\\ 
&=&\int_{S^1} v^*d \iota_{\hat{v}}\zeta+\int_{S^1} v^*\iota_{\hat{v}} d\zeta\\ 
&=&\int_{S^1} d v^*\iota_{\hat{v}}\zeta+\int_{S^1} v^*\iota_{\hat{v}} d\zeta\\ 
&=&\int_{S^1} v^*\iota_{\hat{v}} d\zeta.
\end{eqnarray*}
Alternatively, we can derive this formula using integration by parts as follows
\begin{eqnarray*}
d \mathcal{H}^1(v)\hat{v}&=&\int_{S^1}dF(v)\hat{v}\frac{d}{dt}H(v)dt+\int_{S^1}F(v)\frac{d}{dt}\big(dH(v)\hat{v}\big)dt\\
&=&\int_{S^1}dF(v)\hat{v}\frac{d}{dt}H(v)dt-\int_{S^1}\frac{d}{dt}F(v)dH(v)\hat{v}dt\\
&=&\int_{S^1}dF \wedge dH(v)(\hat{v},\partial_t v)dt.
\end{eqnarray*}
Using this formula we obtain
\begin{eqnarray*}
\big(\mathfrak{a}_\omega-\alpha d\mathcal{H}^1\big)(v)\hat{v}=\int_{S^1} v^*\iota_{\hat{v}}(\omega-\alpha d\zeta).
\end{eqnarray*}
Hence setting
$$\omega^\alpha:=\omega-\alpha d\zeta$$
we can write this as
$$\mathfrak{a}_\omega-\alpha d\mathcal{H}^1=\mathfrak{a}_{\omega^\alpha}.$$
This implies that
\begin{equation}\label{defsym}
\mathfrak{a}_\omega-d\mathcal{H}_{1,\alpha}=\mathfrak{a}_\omega-d\mathcal{H}^0-\alpha d\mathcal{H}^1=
\mathfrak{a}_{\omega^\alpha}-d\mathcal{H}.
\end{equation}
By this formula, the retardation is now incorporated in the deformation of the symplectic form
$\omega$ to the closed two-form $\omega^\alpha$. Note that because nondegeneracy is an open condition, for $\alpha$ small enough the closed two-form $\omega^\alpha$ is still symplectic. 

\subsection{First order in the nonautonomous case and symplectic forms which depend on time}

In the nonautonomous case the discussion of the previous subsection has to be modified in two points. In the first place the
deformation of the symplectic form now becomes a time-dependent two-form. Secondly, the Hamiltonian has to be adjusted as well. 
\\ \\
The formula from the previous subsection gets a bit more involved, namely
\begin{eqnarray*}
\frac{\partial}{\partial \alpha} H_v(\alpha,t)&=&dH_{t+\alpha F(v(t))}\Big(v\big(t+\alpha F(v(t))\big)\Big)\partial_t v
\big(t+\alpha F(v(t))\big)F(v(t))\\
& &+\partial_t H\Big(v\big(t+\alpha F(v(t))\big),t+\alpha F(v(t))\Big)F(v(t))
\end{eqnarray*}
from which we infer
$$\frac{\partial}{\partial \alpha} H_v(0,t)=dH_t(v(t))\partial_t v(t)F(v(t))
+\partial_t H(v(t),t)F(v(t)).$$
We define a smooth family of one-forms on $M$
$$\zeta_t=FdH_t, \quad t \in S^1$$
which give rise to a smooth family of two-forms
$$\omega^\alpha_t:=\omega-\alpha d\zeta_t,\quad t \in S^1.$$
We further define an $\alpha$-dependent periodic Hamiltonian
$$H^\alpha \in C^\infty(M\times S^1,\mathbb{R}), \quad (v,t) \mapsto H(v,t)+\alpha \partial_t H(v,t)F(v).$$
As in (\ref{mean}) we let $\mathcal{H}^\alpha$ be the mean value of $H^\alpha$.
With this notation the analogon of (\ref{defsym}) in the nonautonomous case becomes
\begin{equation}\label{defsym2}
\mathfrak{a}_\omega-d\mathcal{H}_{1,\alpha}=
\mathfrak{a}_{\omega^\alpha}-d\mathcal{H}^\alpha.
\end{equation}
Using Floer homology we shall prove the following theorem about critical points of (\ref{defsym2}). This gives in affirmative answer to Question\,$1^-$ and Question\,$2^-$ in Section~\ref{arnold} in first order in the fine structure constant. 
\begin{thm}\label{main}
Suppose that $(M,\omega)$ is a closed symplectically aspherical manifold. Then there exists $\alpha_0=\alpha_0(H,F)>0$
with the following property. For every $\alpha \in [0,\alpha_0]$ the following statements are true.
\begin{description}
 \item[(i)] All critical points of $\mathfrak{a}_\omega-d\mathcal{H}_{1,\alpha}$ are spacelike.
 \item[(ii)] We have the following lower bound on the number of critical points
 $$\#\mathrm{crit}\big(\mathfrak{a}_\omega-d\mathcal{H}_{1,\alpha}\big) \geq \mathrm{Crit}(M)$$
 where $\mathrm{Crit}(M)$ is the minimal number of critical points of a function on $M$.
 \item[(iii)] Assume that not through every point of $M$ there passes a critical point of 
 $\mathfrak{a}_\omega-d\mathcal{H}_{1,\alpha}$. Then there exist two critical points with different action. 
\end{description}
\end{thm}

\section{Floer homological techniques}

\subsection{Action functional for symplectic forms depending on time}

Assume that $(M,\omega)$ is a connected closed symplectic manifold which is symplectically aspherical, meaning
that the integral of $\omega$ over every sphere vanishes. Suppose that we have a smooth family of one-forms
$\theta_t \in \Omega^1(M)$ for $t \in \mathbb{R}$ which meets the following assumptions.
\begin{description}
 \item[(H1)] The family is twisted periodic in time in the sense that there exists a smooth function
 $K \in C^\infty(M,\mathbb{R})$ satisfying
 $$\theta_{t+1}=\theta_t+dK,\quad \forall\,\,t\in \mathbb{R}.$$
 \item[(H2)] For each $t \in \mathbb{R}$ the two-form
 $$\omega_t:=\omega+d\theta_t \in \Omega^2(M)$$
 is still symplectic. 
\end{description}
Note that it follows from (H1) that 
$$\omega_{t+1}=\omega_t, \quad t \in \mathbb{R}$$
so that we can think of the family of symplectic forms $\omega_t$ as parametrized over the circle. Observe further
that these forms represent the same de Rham cohomology class as $\omega$ and therefore are
themselves symplectically aspherical. Moreover, by Hypothesis (H1) again, the derivative of the family of one-forms $\theta_t$ with respect to the family parameter $t$ is periodic as well, i.e.,
$$\partial_t \theta_{t+1}=\partial_t \theta_t,\quad t\in \mathbb{R}$$
so that we might think of the family of one-forms $\partial_t \theta_t \in \Omega(M)$ as parametrized by the circle as well. 
We interpret the family $\{\theta_t\}_{t\in \mathbb{R}}$ as a one-form on the product manifold
$M \times \mathbb{R}$ which we denote by
$$\theta \in \Omega^1(M\times \mathbb{R}).$$
With respect to the splitting
$$T(M\times \mathbb{R})=TM \times \mathbb{R}$$
its exterior derivative reads at a point $(x,t) \in M \times \mathbb{R}$
$$d\theta_{(x,t)}=(d\theta_t)_x-(\partial_t \theta_t)_x \wedge dt$$
so that we can think of $d \theta$ as a two-form on $M \times S^1$.
If
$$\pi \colon M \times S^1 \to M$$
is the projection along the circle, then we obtain a closed two-form
$$\Omega=\pi^* \omega+d\theta \in \Omega^2(M\times S^1).$$
Abbreviate by
$$\mathcal{L}_0 \subset C^\infty(S^1,M)$$
the component of the free loop space of $M$ consisting of contractible loops. Pick a base point
$p \in M$. Suppose that $v \in \mathcal{L}_0$. Since $v$ is contractible there exists a smooth map
$$\overline{v} \colon [0,1] \times S^1 \to M$$
with the properties that
\begin{equation}\label{filldisk}
\overline{v}(0,t)=p,\,\,t \in S^1,\quad \overline{v}(1,t)=v(t),\,\,t \in S^1.
\end{equation}
We define
$$w_{\overline{v}} \colon [0,1] \times S^1 \to M\times S^1, \quad
(r,t) \mapsto (\overline{v}(r,t),t).$$
Note that
$$\pi \circ w_{\overline{v}}=\overline{v}.$$
\begin{lemma}
The integral $\int_{[0,1]\times S^1} w_{\overline{v}}^*\Omega$ only depends on $v$. 
\end{lemma}
\textbf{Proof:} Using concatenations it suffices to show that
$$\int_{[0,1]\times S^1} w_{\overline{p}}^*\Omega=0$$
for every smooth map
$$\overline{p} \colon [0,1] \times S^1 \to M$$
satisfying
$$\overline{p}(0,t)=\overline{p}(1,t)=p,\quad t \in S^1.$$
Because $\overline{p}$ maps the boundary of the cylinder to the fixed point $p$, we can think of $\overline{p}$ as a map
defined on a sphere. Therefore it follows from the assumption that $(M,\omega)$ is symplectically aspherical that
\begin{equation}\label{asph}
\int \overline{p}^*\omega=0.
\end{equation}
Abbreviate
$$\overline{p}_0 \colon [0,1] \to M,\quad r \mapsto \overline{p}(r,0).$$
Note that
$$\overline{p}_0(0)=p=\overline{p}_0(1)$$
so that we can interpret $\overline{p}_0$ as a map with domain a circle, i.e.
$$\overline{p}_0 \in C^\infty(S^1,M).$$
Using (\ref{asph}), the formula of Stokes, and (H1) we obtain
\begin{eqnarray*}
\int_{[0,1]\times S^1} w_{\overline{p}}^*\Omega&=&\int \overline{p}^* \omega+\int_{[0,1]\times S^1} w_{\overline{p}}^*d\theta\\
&=&\int_{[0,1]\times [0,1]} dw_{\overline{p}}^*\theta\\
&=&-\int_{S^1}\overline{p}_0^*\theta_1+\int_{S^1}\overline{p}_0^*\theta_0\\
&=&-\int_{S^1}\overline{p}_0^* dK\\
&=&-\int_{S^1}d\overline{p}_0^* K\\
&=&0.
\end{eqnarray*}
This finishes the proof of the lemma. \hfill $\square$
\\ \\
In view of the lemma we have a well-defined functional
$$\mathcal{A}_\theta \colon \mathcal{L}_0 \to \mathbb{R},\quad v \mapsto \int_{[0,1]\times S^1} w_{\overline{v}}^*\Omega.$$
To compute its differential, we use the identification
$$T(M\times S^1)=TM \times \mathbb{R}.$$
Note that with respect to this splitting we can write at a point $(x,t) \in M\times S^1$
$$\Omega_{(x,t)}=\omega_x+(d\theta_t)_x-(\partial_t \theta_t)_x \wedge dt=(\omega_t)_x-(\partial_t \theta_t)_x \wedge dt.$$
Hence at a point $v \in \mathcal{L}_0$ the differential of $\mathcal{A}_\theta$ applied to a tangent vector
$\hat{v} \in T_v \mathcal{L}_0$ is
\begin{eqnarray}\label{diffzet}
d\mathcal{A}_\theta(v)\hat{v}&=&\int_0^1\Omega\big((\hat{v},0),(\partial_t v,1)\big)dt\\ \nonumber
&=&\int_0^1 \omega_t(\hat{v},\partial_t v)dt-\int_0^1(\partial_t \theta_t)(\hat{v}(t))dt.
\end{eqnarray}
Because by assumption all members of the family of two-forms $\omega_t$ are symplectic we can define a smooth family
of vector fields
$$Y_t \in \Gamma(TM)$$
implicitly by the requirement
\begin{equation}\label{defy}
\partial_t \theta_t=\omega_t(\cdot,Y_t).
\end{equation}
In view of (\ref{diffzet}) we obtain the following lemma.
\begin{lemma}
Critical points of $\mathcal{A}_\theta$ are in one-to-one correspondence with contractible solutions of the ODE
\begin{equation}\label{critsupham}
\partial_t v(t)=Y_t(v(t)),\quad t \in S^1.
\end{equation}
\end{lemma}
\subsection{Two examples}

\subsubsection{The action functional of classical mechanics}

We discuss how we can recover the action functional of classical mechanics as a special case. Suppose that $H_t \in C^\infty(M,\mathbb{R})$ for $t \in S^1$ is a smooth family of Hamiltonians depending periodically on time. For $t \in \mathbb{R}$ define
\begin{equation}\label{primi}
\overline{H}_t:=\int_0^t H_r dr \in C^\infty(M,\mathbb{R}).
\end{equation}
Then
$$\theta_t:=d\overline{H}_t \in \Omega^1(M)$$
satisfy hypothesis (H1) for
$$K=\overline{H}_1.$$
Moreover, since $\theta_t$ is already exact, we have
$$\omega_t=\omega+d\theta_t=\omega$$
so that hypothesis (H2) is satisfied as well.

We next compute the action functional $\mathcal{A}_\theta$ for this example.
Let $v \in \mathcal{L}_0$ and choose $\overline{v} \colon [0,1]\times S^1 \to M$ satisfying (\ref{filldisk}). Abbreviate
$$\overline{v}_0 \colon [0,1] \to M, \quad r \mapsto \overline{v}(r,0).$$ 
Note that
$$\overline{v}_0(0)=p,\quad \overline{v}_0(1)=v(0)=v(1).$$
Using this notation we obtain using Stokes' theorem
\begin{eqnarray}\label{actex}
\mathcal{A}_\theta(v)&=&\int \overline{v}^*\omega+\int_0^1 \theta_t(\partial_t v)dt+\int_0^1 \theta_0(\partial_t
\overline{v}_0)dt-\int_0^1 \theta_1(\partial_t \overline{v}_0)dt\\ \nonumber
&=&\int \overline{v}^*\omega+\int_0^1 d\overline{H}_t(v(t))\partial_t v dt-\int_0^1 d\overline{H}_1(\partial_t \overline{v}_0)dt\\ \nonumber
&=&\int \overline{v}^*\omega+\int_0^1\frac{d}{dt}\big(\overline{H}_t(v(t))\big)dt-
\int_0^1\big(\partial_t \overline{H}_t\big)(v(t))dt\\ \nonumber
& &-\int_0^1\frac{d}{dt}\overline{H}_1(\overline{v}_0(t))dt\\ \nonumber
&=&\int \overline{v}^*\omega+\overline{H}_1(v(1))-\overline{H}_0(v(0))-\int_0^1 H_t(v(t))dt-\overline{H}_1(v(1))\\ \nonumber
& &+\overline{H}_1(p)\\ \nonumber
&=&\int \overline{v}^*\omega-\int_0^1 H_t(v(t))dt+\overline{H}_1(p)\\ \nonumber
&=:&\mathcal{A}_H(v)+\overline{H}_1(p).
\end{eqnarray}
The result is just the action functional of classical mechanics up to addition of the constant $\overline{H}_1(p)$.
Note that addition of a constant has no influence on the critical points and therefore the critical points of
$\mathcal{A}_\theta$ in this example are the critical points of the action functional of classical mechanics, which are the periodic orbits of the Hamiltonian vector field of $H$, i.e., solutions of the ODE
\begin{equation}\label{critham}
\partial_t v(t)=X_{H_t}(v(t)), \quad t \in S^1.
\end{equation}
This can be as well verified directly: In this example
$$\partial_t \theta_t=d\partial_t \overline{H}_t=dH_t$$
so that from (\ref{defy}) we infer that
$$Y_t=X_{H_t}.$$
It now follows from (\ref{critsupham}) that the critical points of $\mathcal{A}_\theta$ are the solutions of (\ref{critham}).

\subsubsection{First order Neumann one-forms}\label{firstorder}

We discuss how the Taylor polynomial of a Neumann one-form in first order in the fine structure constant on a symplectically aspherical symplectic manifold for small enough fine structure constant arises as the differential of a functional $\mathcal{A}_ \theta$. For that purpose suppose that we have a smooth loop of one-forms
$$\zeta_t \in \Omega^1(M),\quad t \in S^1$$
with the property that
$$\omega_t=\omega+d\zeta_t \in \Omega^2(M)$$
are symplectic for every $t \in S^1$. Note that different to the case of a one parameter family of one-forms $\theta_t$ for which we required in
(H1) just twisted periodicity we require for the one parameter family of one-forms $\zeta_t$ periodicity in $t$. For 
a Hamiltonian $H \in C^\infty(M\times S^1,\mathbb{R})$ we define $\overline{H} \in C^\infty(M \times \mathbb{R},\mathbb{R})$ by (\ref{primi}) and set
$$\theta_t=\zeta_t+d\overline{H}_t.$$
Abbreviate
$$\mathcal{A}_{\zeta,H}\colon \mathcal{L}_0 \to \mathbb{R},\quad \mathcal{A}_{\zeta,H}:=\mathcal{A}_\zeta-\mathcal{H}$$
where $\mathcal{H}$ is the mean value of $H$ as in (\ref{mean}). Then the same computation as in (\ref{actex}) gives the formula
$$\mathcal{A}_\theta=\mathcal{A}_{\zeta,H}+\overline{H}_1(p).$$
From this expression we learn that $\mathcal{A}_\theta$ and $\mathcal{A}_{\zeta,H}$ have the same critical points. 

\subsection{Gradient flow lines}

Suppose that $J_t$ is a smooth family of almost complex structures on $M$ such that $J_t$ is $\omega_t$ compatible, meaning that
$$g_t(\cdot,\cdot)=\omega_t(\cdot,J_t \cdot)$$
is a Riemannian metric on $M$. That such families exist is for example explained in \cite[Chapter 4]{mcduff-salamon1}. Integrating $g_t$ gives rise to an $L^2$-metric $g$ on the free loop space of $M$. From
(\ref{diffzet}) and (\ref{defy}) we infer that the gradient of $\mathcal{A}_\theta$ with respect to $g$ at a point
$v \in \mathcal{L}_0$ is
$$\nabla \mathcal{A}_\theta(v)=J_t(v)(\partial_t v -Y_t(v)).$$
Hence negative gradient flow lines of $\nabla \mathcal{A}_\theta$, which are formally maps
$$v \colon \mathbb{R} \to \mathcal{L}_0$$
satisfying
\begin{equation}\label{abgrad}
\partial_s v(s)+\nabla \mathcal{A}_\theta(v(s))=0,\quad s \in \mathbb{R},
\end{equation}
can be interpreted as maps
$v \in C^\infty(\mathbb{R} \times S^1,M)$ satisfying the following perturbed nonlinear Cauchy Riemann equation
\begin{equation}\label{grad}
\partial_s v(s,t)+J_t(v(s,t))\big(\partial_t v(s,t)-Y_t(v(s,t))\big),\quad (s,t) \in \mathbb{R} \times S^1.
\end{equation}
The energy of a solution of (\ref{grad}) is defined as
$$E(v)=\int_{-\infty}^\infty g(\partial_s v,\partial_s v)ds=\int_{-\infty}^\infty \int_0^1 g_t(\partial_s v(s,t),\partial_s v(s,t))dsdt.$$
Because for every $t \in S^1$ the symplectic form $\omega_t$ is symplectically aspherical there is no bubbling and the usual compactness result in Floer homology, see for instance \cite{audin-damian, mcduff-salamon2}, implies the following theorem.
\begin{thm}\label{gromcom}
Suppose that $v_\nu$ for $\nu \in \mathbb{N}$ is a sequence of gradient flow lines, i.e., solutions of (\ref{grad}), whose energy is uniformly bounded in the sense that there exists a constant $c$ such that
$$E(v_\nu) \leq c,\quad \nu \in \mathbb{N}.$$
Then there exists a subsequence $\nu_j$ and another gradient flow line $v$ of energy $E(v)\leq c$ such that the sequence
$v_{\nu_j}$ converges in the $C^\infty_{\mathrm{loc}}$-topology to $v$.
\end{thm}
We remark that in general one only has local convergence and not global convergence due to the phenomenon that gradient flow lines might break. 

\subsection{Stretching homotopies}

In this subsection we assume that the family of one-forms $\theta_t$ satisfies besides hypothesis (H1) the following 
sharpened version of hypothesis (H2). 
\begin{description}
 \item[(H2+)] For each $t \in S^1$ and $\tau \in [0,1]$ the two-form
 $$\omega_{\tau,t}:=\omega+\tau d\theta_t \in \Omega^2(M)$$
 is symplectic. 
\end{description}
Suppose further that $\beta \in C^\infty(\mathbb{R},[0,1])$ is a bump function satisfying
$$\beta(s)=1,\quad s \in [-1,1],\qquad \beta(s)=0,\quad |s|\geq 2.$$
We choose a smooth family of bump functions
$$\beta_r \in C^\infty(\mathbb{R},[0,1]), \quad r \in [0,\infty)$$
with the property that $\beta_0$ vanished identically, i.e.
$$\beta_0(s)=0, \quad s \in \mathbb{R}$$
and for $r \geq 1$ the bump functions $\beta_r$ are obtained by stretching the fixed bump function $\beta$ in the sense that
$$\beta_r(s)=\left\{\begin{array}{cc}
\beta(s+r) & s \leq -r-1\\
1 & -r-1 \leq s \leq r+1\\
\beta(s-r) & s \geq 1+r.
\end{array}
\right.$$
We now define a family of time-dependent action functionals. Namely for $r \in [0,\infty)$ we define
$$\mathcal{A}_r \colon \mathcal{L}_0 \times \mathbb{R} \to \mathbb{R}, \quad
(v,s) \mapsto \mathcal{A}_{\beta_r(s)\theta}(v).$$
Note that all the  action functionals are asympotically independent of the variable $s$ and given by the unperturbed area functional $\mathcal{A}_0$. However, as $r$ increases the action functionals $\mathcal{A}_r$ coincide on an increasing part of the real line with the action functional $\mathcal{A}_\theta$. The family $\mathcal{A}_r$ therefore interpolates between the area functional $\mathcal{A}_0$ and $\mathcal{A}_\theta$.

Since by assumption (H2+) all two-forms $\omega_{\tau,t}$ for $t \in S^1$ and $\tau \in [0,1]$  are symplectic, we can choose a smooth family $J_{\tau,t}$ of $\omega_{\tau,t}$-compatible almost complex structures. This gives rise to a smooth two-parameter family of Riemannian metrics on $M$
$$g_{\tau,t}(\cdot,\cdot)=\omega_{\tau,t}(\cdot,J_{\tau,t}\cdot).$$
Integrating with respect to the $t$-variable gives rise to a one-parameter family of Riemannian metrics $g_\tau$ on the free loop space $\mathcal{L}$ and we denote by
$$\nabla_\tau=\nabla_{g_\tau}$$
the gradient with respect to the metric $g_\tau$. We consider the following problem. We are looking for pairs
$(v,r)$ where $r \in [0,\infty)$ and $v \colon \mathbb{R} \to \mathcal{L}_0$ which satisfy
\begin{equation}\label{homhom}
\partial_s v(s)+\nabla_{\beta_r(s)}\mathcal{A}_r\big(v(s),s\big)=0,\,\,s \in \mathbb{R},\qquad
E(v,r)<\infty.
\end{equation}
Here the energy is defined as
$$E(v,r)=\int_{-\infty}^\infty g_{\beta_r(s)}\big(\partial_s v(s),\partial_s v(s)\big)ds.$$
If one interprets solutions of the $s$-dependent gradient flow equation in (\ref{homhom}) as in (\ref{grad}) as maps from the cylinder to $M$, they are solutions of the PDE
\begin{equation}\label{homgrad}
\partial_s v(s,t)+J_{\beta_r(s),t}(v(s,t))\big(\partial_t v(s,t)-Y_{\beta_r(s),t}(v(s,t))\big)=0,\quad (s,t) \in \mathbb{R} \times S^1,
\end{equation}
where for $t \in S^1$ and $\tau \in [0,1]$ the vector field $Y_{\tau,t}$ is implicitly defined by the equation
$$\tau \partial_t \theta_t=\omega_{\tau,t}(\cdot,Y_{\tau,t}).$$
Remember that all the time dependent action functionals $\mathcal{A}_{\beta_r(s)}$ are asymptotically equal to the unperturbed area functional $\mathcal{A}_0$. The only critical points of the area functional are constant loops. Indeed, the area
functional is Morse-Bott with critical set diffeomorphic to $M$, by interpreting constant loops as points in $M$. Therefore
solutions of (\ref{homhom}) exponentially converge to constants. We abbreviate
\begin{equation}\label{moduli1}
\mathcal{M}=\big\{(v,r)\,\,\textrm{solutions of (\ref{homhom})}\big\}
\end{equation}
the moduli space of solutions of the problem (\ref{homhom}).
\begin{prop}\label{uniprop}
There exists a finite constant $c$ which bounds the energy of all members of the moduli space $\mathcal{M}$ uniformly, i.e.,
$$E(v,r) \leq c, \quad (v,r) \in \mathcal{M}.$$
\end{prop}
\textbf{Proof: } Suppose that $(v,r) \in \mathcal{M}$. Using the gradient flow equation we rewrite the energy as follows
\begin{eqnarray}\label{u1}
E(v,r)&=&\int_{-\infty}^\infty g_{\beta_r(s)}\big(\partial_s(v(s)),\partial_s v(s)\big)ds\\ \nonumber
&=&-\int_{-\infty}^\infty g_{\beta_r(s)}\Big(\nabla_{\beta_r(s)}\mathcal{A}_{\beta_r(s)\theta}(v(s)),\partial_s v(s)\Big)ds\\
\nonumber
&=&-\int_{-\infty}^\infty d\mathcal{A}_{\beta_r(s) \theta}(v(s))\partial_s v(s) ds\\ \nonumber
&=&-\int_{-\infty}^\infty \frac{d}{ds}\mathcal{A}_{\beta_r(s) \theta}(v(s))ds+\int_{-\infty}^\infty
\big(\partial_s \mathcal{A}_{\beta_r(s)\theta}\big)(v(s))ds\\ \nonumber
&=&\int_{-\infty}^\infty
\big(\partial_s \mathcal{A}_{\beta_r(s)\theta}\big)(v(s))ds.
\end{eqnarray}
For the last equation we used the fact that by the assumption that $E(v,r)$ is finite $v(s)$ converges asymptotically to constant loops at which the area functional vanishes. To continue with the proof we need the following lemma. All Riemannian metrics on $M$ are equivalent. Therefore in the statement of the following lemma it does not matter what metric we use 
to define the $L^1$-norm. However, in order that we do not need to alter the constant later on, we take the metric
$g_{\beta_r(s)}$ on the free loop space. 
\begin{lemma}\label{unilem}
There exists a finite constant $c_1$ such that for $v \in \mathcal{L}_0$, $s \in \mathbb{R}$, and $r \in [0,\infty)$ 
the following estimate holds
$$\big|\big(\partial_s \mathcal{A}_{\beta_r(s)\theta}\big)(v)\big| \leq c_1\big(||\partial_t v||_{L_1}+1\big).$$
\end{lemma}
\textbf{Proof: } Choose a smooth function $\overline{v} \colon [0,1] \times S^1 \to M$ satisfying (\ref{filldisk}).
Abbreviate
$$\overline{v}_0 \colon [0,1] \to M, \quad r \mapsto \overline{v}(r,0).$$
Note that $\overline{v}_0$ is a smooth path from the basepoint $p \in M$ to $v(0)$ the starting and ending point of the 
loop $v$. We then have
\begin{eqnarray*}
\big(\partial_s \mathcal{A}_{\beta_r(s)\theta}\big)(v)&=&\beta'_r(s)\bigg(
\int_0^1 \theta_t(\partial_t v)dt+\int_0^1\theta_0(\partial_t \overline{v}_0)dt-\int_0^1 \theta_1(\partial_t \overline{v}_0)dt\bigg)\\ 
&=&\beta'_r(s)\bigg(
\int_0^1 \theta_t(\partial_t v)dt-\int_0^1 dK(\partial_t \overline{v}_0)dt\bigg)\\
&=&\beta'_r(s)\bigg(
\int_0^1 \theta_t(\partial_t v)dt-\int_0^1 \frac{d}{dt}K(\overline{v}_0)dt\bigg)\\
&=&\beta'_r(s)\bigg(
\int_0^1 \theta_t(\partial_t v)dt+K(p)-K(v(0))\bigg).
\end{eqnarray*}
Note that we have chosen the family $\beta_r$ which depends on the parameter $r$ living in the noncompact space 
$[0,\infty)$ in such a way that the derivative $\beta'_r(s)$ is uniformly bounded. Therefore the assertion of the lemma follows from the formula above. \hfill $\square$
\\ \\
\textbf{Proof of Proposition~\ref{uniprop} continued: } If $(v,r) \in \mathcal{M}$, then $v$ is a solution of the PDE (\ref{homgrad}) from which we infer that
$$\partial_t v(s,t)=J_{\beta_r(s),t}(v(s,t))\partial_s v(s,t)+Y_{\beta_r(s),t}(v(s,t)).$$
Combining this with Lemma~\ref{unilem} and using that $Y_{\beta_r(s),t}$ is uniformly bounded in $M$, $r$, $s$, and $t$ we obtain a constant $c_2$ with the property that 
\begin{equation}\label{u2}
\big|\big(\partial_s \mathcal{A}_{\beta_r(s)\theta}\big)(v(s))\big| \leq c_2\big(||\partial_s v(s)||_{L_1}+1\big).
\end{equation}
For $r \in [0,\infty)$ abbreviate
$$C_r:=\mathrm{supp} \beta'_r.$$
By our construction of the cut-off functions $\beta_r$ there exists a constant $c_3$ independent of
$r$ such that
\begin{equation}\label{u3}
|C_r| \leq c_3
\end{equation}
where $|C_r|$ denotes the Lebesgue measure of the set $C_r \subset \mathbb{R}$.
Using (\ref{u1}),(\ref{u2}), and (\ref{u3}) we estimate
\begin{eqnarray*}
E(v,r)&\leq& \int_{-\infty}^\infty\big|\big(\partial_s \mathcal{A}_{\beta_r(s)\theta}\big)(v(s))\big|ds\\
&=&\int_{C_r}\big|\big(\partial_s \mathcal{A}_{\beta_r(s)\theta}\big)(v(s))\big|ds\\
&\leq&c_2 \int_{C_r}\big(||\partial_s v(s)||_{L_1}+1\big)ds\\
&\leq&c_2 \big|\big|\partial_s v|_{C_r \times S^1}\big|\big|_{L_1}+c_2 c_3\\
&\leq& c_2\sqrt{\big|C_r\times S^1\big|}\,\big|\big|\partial_s v|_{C_r \times S^1}\big|\big|_{L_2}+c_2 c_3\\
&\leq& c_2\sqrt{c_3}\sqrt{||\partial_s v||^2_{L_2}}+c_2c_3\\
&=&c_2\sqrt{c_3}\sqrt{E(v,r)}+c_2c_3.
\end{eqnarray*}
Abbreviating
$$c_4:=c_2 \cdot \max\{c_3,\sqrt{c_3}\}$$
we deduce from this inequality the more condensed inequality
\begin{equation}\label{u4}
E(v,r) \leq c_4\big(\sqrt{E(v,r)}+1\big).
\end{equation}
We claim that (\ref{u4}) implies that 
\begin{equation}\label{u5}
E(v,r) \leq \max\{4c_4^2,1\}.
\end{equation}
To see that we can assume without loss of generality that
$$E(v,r) \geq 1$$
and therefore as well
$$\sqrt{E(v,r)}\geq 1.$$
Plugging this inequality into (\ref{u4}) we obtain the inequality
$$E(v,r) \leq 2c_4 \sqrt{E(v,r)}$$
or dividing both sides by $\sqrt{E(v,r)}$
$$\sqrt{E(v,r} \leq 2c_4$$
from which (\ref{u5}) follows by squaring. By setting
$$c=\max\{4c_4^2,1\}$$
the assertion of the proposition follows now from (\ref{u5}). \hfill $\square$

\subsection{Uniruledness}
In this subsection we study the moduli space of finite energy gradient flow lines of $\nabla \mathcal{A}_\theta$, i.e.,
$$\mathcal{N}:=\big\{v\,\,\textrm{solution of (\ref{abgrad}}), E(v)<\infty\big\}.$$ 
Note that each critical point of $\mathcal{A}_\theta$ interpreted as a constant gradient flow line belongs to the moduli space $\mathcal{N}$ as well. Indeed, a gradient flow line corresponds to a critical point if and only if its energy is zero. For $q \in M$ we abbreviate
$$\mathcal{N}_q:=\big\{v \in \mathcal{N}, v(0,0)=q\big\}$$
the moduli space of all finite energy gradient flow lines which pass at $(0,0)$ through $q$. In this subsection we prove the following uniruledness result.
\begin{thm}\label{uniruled}
Under hypotheses (H1) and (H2+) for every $q \in M$ the moduli space $\mathcal{N}_q$ is not empty. 
\end{thm}
The proof of Theorem~\ref{uniruled} uses heavily the moduli space $\mathcal{M}$ introduced in (\ref{moduli1}). We have a natural projection
$$\pi \colon \mathcal{M} \to [0,\infty), \quad (v,r) \mapsto r$$
and we abbreviate
$$\mathcal{M}_r:=\pi^{-1}(r)$$
so that we have
$$\mathcal{M}=\bigcup_{r \in [0,\infty)}\mathcal{M}_r.$$
We might think that the moduli space $\mathcal{N}$ contains a kind of boundary of the moduli space $\mathcal{M}$ as
$r$ goes to infinity. Therefore in order to proof Theorem~\ref{uniruled} we need to find suitable elements in the moduli space $\mathcal{M}_r$ for arbitrarily large parameter $r$. Here is how this works.
\\ \\
\textbf{Proof of Theorem~\ref{uniruled}: }For $q \in M$ we abbreviate
$$\mathcal{M}_q:=\big\{(v,r) \in \mathcal{M}, v(0,0)=q\big\}.$$
Set
$$\pi_q \colon \mathcal{M}_q \to [0,\infty),\quad (v,r) \mapsto r$$
and define for $r$ and $R$ nonnegative real numbers
$$\mathcal{M}_{q,r}:=\pi_q^{-1}(r), \quad \mathcal{M}_q^R:=\pi_q^{-1}\big([0,R]\big).$$
Note that
$$\mathcal{M}_q^R=\bigcup_{r \in [0,R]}\mathcal{M}_{q,r}.$$
To prove the theorem we need the following proposition.
\begin{prop}\label{nonempty}
For each $r \in [0,\infty)$ the moduli space $\mathcal{M}_{r,q}$ is nonempty.
\end{prop}
The proof of this proposition hinges on the following proposition.
\begin{prop}\label{compact}
For each $R \in [0,\infty)$ the moduli space $\mathcal{M}_q^R$ is compact. 
\end{prop}
\textbf{Proof: } By Proposition~\ref{uniprop} the energy of all gradient flow lines in $\mathcal{M}_q^R$ is uniformly bounded. Hence by Theorem~\ref{gromcom} each sequence in $\mathcal{M}_q^R$ has a locally converging subsequence. To see that it converges globally note that all action functionals $\mathcal{A}_r$ are asymptotically equal to the area functional whose only critical points are constants at which the action vanishes. Therefore gradient flow lines cannot break and the convergence is global. This finishes the proof of the proposition. \hfill $\square$
\\ \\
\textbf{Proof of Proposition~\ref{nonempty}: }For $R \in [0,\infty)$ the moduli space $\mathcal{M}_q^R$ can be interpreted as the zero set of a Fredholm section of Fredholm index one from a Hilbert manifold with boundary to a Hilbert bundle over it. If this section is transverse then $\mathcal{M}_q^R$ is a one dimensional manifold with boundary. The boundary of this manifold is
$$\partial \mathcal{M}_q^R=\mathcal{M}_{q,0} \cup \mathcal{M}_{q,R}.$$
The moduli space $\mathcal{M}_{q,0}$ consists of a single element, namely the constant gradient flow line $q$. By Proposition~\ref{compact} the moduli space $\mathcal{M}_q^R$ is compact as well. A compact one-dimensional manifold with boundary consists of a finite union of circles and closed intervals. In particular, the cardinality of its boundary points is even. Therefore the moduli space $\mathcal{M}_{q,R}$ cannot be empty. This proves the proposition in the case the Fredholm section is transverse. If it is not transverse we slightly perturb it to make it transverse and the proposition follows by the same argument again. \hfill $\square$
\\ \\
\textbf{Proof of Theorem~\ref{uniruled} continued: } Choose a sequence $r_\nu$ of nonnegative real numbers converging to infinity as $\nu$ goes to infinity. By Proposition~\ref{nonempty} there exists a gradient flow line $v_\nu$ of 
$\mathcal{A}_{r_\nu}$ of finite energy through $q$ for every $\nu \in \mathbb{N}$. By Proposition~\ref{uniprop} the energy of the sequence of gradient flow lines $v_\nu$ is uniformly bounded. Therefore by Theorem~\ref{gromcom} it admits a subsequence which locally converges to a gradient flow line of $\mathcal{A}_\theta$. In particular, the moduli space 
$\mathcal{N}_q$ is not empty. This finishes the proof of the theorem. \hfill $\square$
\subsection{Existence of two critical points of different action and proof of the Arnold conjecture}

In this subsection we prove the following theorem.
\begin{thm}\label{zwei}
Assume (H1) and (H2+) and suppose further that not through every point $q \in M$ passes a critical point of
$\mathcal{A}_\theta$. Then there exist two critical points of $\mathcal{A}_\theta$ of different action. 
\end{thm}
\textbf{Proof: } By the assumption of the theorem we can choose $q \in M$ such that no critical point of $\mathcal{A}_\theta$ goes through $q$. By Theorem~\ref{uniruled} there exists a gradient flow line $w$ of $\mathcal{A}_\theta$ of finite energy through $q$. By our choice of $q$ the gradient flow line $w$ cannot be a critical point so that its energy is positive, i.e.
$$E(w) \in (0,\infty).$$
Abbreviate by
$$\Omega^-(w),\,\,\Omega^+(w) \subset \mathrm{crit}(\mathcal{A}_\theta)$$
the Omega limit sets of $w$. Namely $\Omega^+(w)$ consists of all critical points $v$ of $\mathcal{A}_\theta$ for which there exists a sequence $s_\nu$ converging to $\infty$ such that 
$$v=\lim_{\nu \to \infty} w(s_\nu)$$
and similarly for $\Omega^-(w)$ if the sequence instead converges to $-\infty$. Because the energy of $w$ is finite we infer from Theorem~\ref{gromcom} that 
$$\Omega^+(w) \neq \emptyset,\quad \Omega^-(w) \neq \emptyset.$$
In particular, choose
$$v^+ \in \Omega^+(w),\quad v^- \in \Omega^-(w).$$ 
Using the gradient flow equation (\ref{abgrad}) we obtain
\begin{eqnarray*}
0 &<& E(w)\\
&=&\int_{-\infty}^\infty g(\partial_s w,\partial_s w)ds\\
&=&-\int_{-\infty}^\infty g\big(\nabla \mathcal{A}_\theta(w),\partial_s w\big)ds\\
&=&-\int_{-\infty}^\infty d\mathcal{A}_\theta(w)\partial_s w ds\\
&=&-\int_{-\infty}^\infty \frac{d}{ds}\mathcal{A}_\theta(w) ds\\
&=&\mathcal{A}_\theta(v^-)-\mathcal{A}_\theta(v^+).
\end{eqnarray*}
In particular,
$$\mathcal{A}_\theta(v^-) \neq \mathcal{A}_\theta(v^+).$$
Hence we found two critical points of $\mathcal{A}_\theta$ of different action and the theorem is proven. \hfill $\square$
\\ \\
The results established to prove existence of two critical points of different action can also be used to give a proof of the original Arnold conjecture in this set-up.
\begin{thm}\label{arth}
Assume hypotheses (H1) and (H2+). Then 
$$\#\mathrm{crit}(\mathcal{A}_\theta) \geq \mathrm{Crit}(M)$$
where $\mathrm{Crit}(M)$ denotes the minimal number of critical points of a function on $M$.
\end{thm}
\textbf{Proof: }The results to prove Theorem~\ref{uniruled} combined with the techniques in \cite{hofer} show that the
evaluation map 
$$\mathrm{ev} \colon \mathcal{N} \to M, \quad v \mapsto v(0,0)$$
induces an injection
$$\mathrm{ev}^* \colon H^*(M;\mathbb{Z}_2) \to H^*(\mathcal{N};\mathbb{Z}_2).$$
Taking coherent orientations into account as established in \cite{floer-hofer}, this statement can be improved to get an injection
$$\mathrm{ev}^* \colon H^*(M;G) \to H^*(\mathcal{N};G)$$
for every coefficient group $G$. Since $M$ is symplectically aspherical the results in \cite{rudyak, rudyak-oprea} imply the theorem. \hfill $\square$

\subsection{Proof of Theorem~\ref{main}}

Combining formula (\ref{defsym2}) with the discussion of Paragraph~\ref{firstorder} we see that the Floer theoretic set-up developed in this section can be applied to the detection of critical points appearing in Theorem~\ref{main}. Moreover, for small enough fine structure constant hypothesis (H2+) is satisfied. From (\ref{critsupham}) we infer that 
for small enough fine structure constant we can assume in addition that all critical points are spacelike and therefore assertion (i) of the theorem holds true. Assertion (ii) follows from Theorem~\ref{arth} and assertion (iii) follows from
Theorem~\ref{zwei}. This finishes the proof of the theorem. \hfill $\square$

\appendix

\section{Fuzzy Neumann one-forms}

Fuzzy Neumann one-forms arise when the retardation is not local anymore but is allowed to depend on the whole trajectory. 

\subsection{Motivation}

We were able to prove an Arnold-type conjecture linearly in the fine structure constant. The crucial point here was the observation that the linear term gives rise to a time-dependent perturbation of the symplectic form. If one studies 
quadratic Taylor polynomials in the fine structure constant then additionally terms quadratic in the velocity arise. These terms are of higher order than the symplectic terms and if it is possible to deal with them then complicated analytical tools have to be developed. A possible tool is to allow that the retardation depends as well on the trajectory. 
In case it is possible to prove existence of critical points in this set-up one then could study a limit case 
where in the limit the retardation just depends pointwise on the trajectory. 
\\ \\
The problem with quadratic terms in the velocity is that if the velocity has a peak at some time, taking square makes the peak even much worse. On the other hand if the peak is spread out this effect is much more harmless as the following straightforward estimate shows. For that let $\beta$ be a real valued function on the circle.  

\begin{eqnarray*}
\int_0^1\bigg|\partial_t v(t)\int_0^1 \partial_t v(t+\tau)\beta(\tau)d\tau\bigg|dt&=&
\int_0^1\bigg|\partial_t v(t)\int_0^1 \partial_\tau v(t+\tau)\beta(\tau)d\tau\bigg|dt\\
&=&\int_0^1\bigg|\partial_t v(t)\int_0^1 v(t+\tau)\partial_\tau\beta(\tau)d\tau\bigg|dt\\
&\leq&\int_0^1 |\partial_t v(t)|\int_0^1 |v(t+\tau)| \cdot |\partial_\tau\beta(\tau)|d\tau dt\\
&\leq&||v||_\infty\int_0^1 |\partial_t v(t)|\int_0^1 |\partial_\tau\beta(\tau)|d\tau dt\\
&=&||v||_\infty\cdot ||\partial_t v||_1 \cdot ||\partial_t \beta||_1.
\end{eqnarray*}

\subsection{Fuzzy Neumann one-forms and its first order approximation}

In this subsection we introduce fuzzy Neumann one-forms and its Taylor polynomials and we study the linear Taylor polynomial of a fuzzy Neumann one-form. Given a symplectic form on a manifold we also obtain a symplectic form on its loop space by integrating it. If the symplectic form on the loop space arises in this way then we refer to it as local. We see that the perturbation in first order in the fine structure constant of a fuzzy Neumann one-form can be interpreted as a nonlocal perturbation of the symplectic form on the free loop space. To make the formulas not too complicated we just discuss the autonomous case. Hence we assume that we have a time independent Hamiltonian $H \in C^\infty(M,\mathbb{R})$ and
a retardation $F \in C^\infty(M,\mathbb{R})$. As an additional data we need
$$\beta \in C^\infty(S^1,\mathbb{R}).$$
We introduce
$$\mathcal{H}_{\alpha,\beta}(v):=\int_0^1\int_0^1 H\Big(v\big(t+\alpha F(v(t+\tau))\big)\Big)\beta(\tau)dt d\tau.$$
We abbreviate
$$H_{v,\beta} \in C^\infty([0,\infty)\times S^1,\mathbb{R}),\quad (\alpha,t) \mapsto \int_0^1
H\Big(v\big(t+\alpha F(v(t+\tau))\big)\Big)\beta(\tau)d\tau$$
For $k \in \mathbb{N}_0$ we define
$$\mathcal{H}^k_\beta(v):=\frac{1}{k!}\int_0^1 \frac{\partial^k}{\partial \alpha^k}H_{v,\beta}(0,t)dt$$ 
and for $n \in \mathbb{N}_0$ we set
$$\mathcal{H}_{n,\alpha,\beta}:=\sum_{k=0}^n \alpha^k \mathcal{H}^k_\beta \colon \mathcal{L}\to \mathbb{R}.$$
A \emph{fuzzy Neumann one-form} is then the one-form on the free loop space given by
$$\mathfrak{a}_\omega-d\mathcal{H}_{\alpha,\beta} \in \Omega^1(\mathcal{L})$$
and its Taylor polynomial of degree $n$ with respect to the fine structure constant is
$$\mathfrak{a}_\omega-d\mathcal{H}_{n,\alpha,\beta} \in \Omega^1(\mathcal{L}).$$
Note that if $\beta$ is the delta-distribution surging at zero we obtain the usual Neumann one-forms back. 
\\ \\
To describe the Taylor polynomial of degree one we compute
\begin{eqnarray*}
& &\frac{\partial}{\partial \alpha} H_{v,\beta}(\alpha,t)\\
&=&\int_0^1 dH\Big(v\big(t+\alpha F(v(t+\tau))\big)\Big)\partial_t v
\big(t+\alpha F(v(t+\tau))\big)F(v(t+\tau))\beta(\tau)d\tau.
\end{eqnarray*}
Evaluating this at $\alpha=0$ we obtain the expression
\begin{eqnarray*}
\frac{\partial}{\partial \alpha} H_{v,\beta}(0,t)&=&\int_0^1 F(v(t+\tau))dH(v(t))\partial_t v(t)\beta(\tau)d\tau\\
&=&\int_0^1 F(v(t+\tau))\frac{d}{dt}H(v(t))\beta(\tau)d\tau.
\end{eqnarray*}
Hence 
$$\mathcal{H}_\beta^1(v)=\int_0^1\int_0^1 F(v(t+\tau))\frac{d}{dt}H(v(t))\beta(\tau)d\tau dt$$
and its differential computes to be
\begin{eqnarray*}
d \mathcal{H}^1_\beta(v)\hat{v}&=&\int_0^1 \int_0^1dF(v(t+\tau))\hat{v}(t+\tau)\frac{d}{dt}\big(H(v(t))\big)\beta(\tau) d\tau dt\\
& &+\int_0^1\int_0^1 F(v(t+\tau))\frac{d}{dt}\big(dH(v(t))\hat{v}(t)\big)\beta(\tau) d\tau dt\\
&=&\int_0^1 \int_0^1dF(v(t))\hat{v}(t)\frac{d}{dt}\big(H(v(t-\tau))\big)\beta(\tau) d\tau dt\\
& &-\int_0^1\int_0^1 \frac{d}{dt}\big(F(v(t+\tau))\big)dH(v(t))\hat{v}(t)\beta(\tau) d\tau dt\\
&=&\int_0^1 \int_0^1dF(v(t))\hat{v}(t)dH(v(t-\tau))\partial_t v(t-\tau)\beta(\tau) d\tau dt\\
& &-\int_0^1\int_0^1 dF(v(t+\tau))\partial_t v(t+\tau)dH(v(t))\hat{v}(t)\beta(\tau) d\tau dt.
\end{eqnarray*}
By integrating the symplectic form $\omega$ on $M$ we obtain a symplectic form on the loop space
$$\varpi \in \Omega^2(\mathcal{L})$$
which is given for two tangent vectors $\hat{v}_1, \hat{v}_2 \in T_v \mathcal{L}$ at a point $v \in \mathcal{L}$ by
$$\varpi(\hat{v}_1,\hat{v}_2):=\int_0^1 \omega\big(\hat{v}_1(t),\hat{v}_2(t)\big)dt.$$
We further introduce the nonlocal one-form
$$\varsigma_\beta \in \Omega^1(\mathcal{L})$$
which for a tangent vector $\hat{v} \in T_v \mathcal{L}$ at a point $v \in \mathcal{L}$ is defined as
$$\varsigma_\beta(\hat{v})=\int_0^1 \int_0^1 F(v(t+\tau))dH(v(t))\hat{v}(t)\beta(\tau)d\tau dt.$$
For $\alpha \in [0,\infty)$ we then obtain a family of closed nonlocal two-forms
$$\varpi^\alpha_\beta \in \Omega^2(\mathcal{L}), \quad \varpi^\alpha_\beta:=\varpi-\alpha d\varsigma_\beta.$$
Note that for $\alpha$ small enough this two-form is a nonlocal symplectic form on the free loop space. 
We define the family of closed one-forms
$$\mathfrak{a}_{\omega,\alpha,\beta} \in \Omega^1(\mathcal{L})$$
which at $v \in \mathcal{L}$ are given for $\hat{v} \in T_v \mathcal{L}$ by
$$\mathfrak{a}_{\omega,\alpha,\beta}(v)\hat{v}=\varpi^\alpha_\beta(\hat{v},\partial_t v).$$
With this notation the analogon of formula (\ref{defsym}) in the fuzzy case becomes
$$\mathfrak{a}_\omega-d\mathcal{H}_{1,\alpha,\beta}=\mathfrak{a}_{\omega,\alpha,\beta}-d\mathcal{H}_{0,\beta}.$$
Note that in contrast to (\ref{defsym}) the closed two-from on the free loop space is now not local anymore. To the author's knowledge nobody constructed a Floer homology in this case. However, the author expects that an analogue of Theorem~\ref{main} also holds true in this case and hopes that establishing such a theorem gives useful insight how a nonlocal Floer homology as propagated in \cite{albers-frauenfelder-schlenk1,albers-frauenfelder-schlenk2,albers-frauenfelder-schlenk3} looks like. Using ideas from the previous subsection one can try to see if such a Floer homology can as well be used to tell something about the critical points of Taylor polynomials of Neumann one-forms of higher degree.

\end{document}